\documentclass[11pt,draft]{article}
\usepackage[english]{babel}
\usepackage{amssymb,amsfonts}
\def\BBox{\vrule height 0.5em width 0.6em depth 0em}
\newtheorem{theo}{Theorem}[section]
\newtheorem{lema}{Lemma}[section]
\newtheorem{prop}{Proposition}[section]

\newcommand{\refe}[1]{(\ref{#1})}

\newcommand{\tsig}{\tilde{\sigma}}
\newcommand{\ttau}{\tilde{\tau}}
\newcommand{\bq}{\begin{equation}}
\newcommand{\eq}{\end{equation}}
\newcommand{\dst}{\displaystyle}
\newcommand{\ba}{\begin{array}}
\newcommand{\ea}{\end{array}}

\title{Factorization method for difference
equations of hypergeometric type on nonuniform lattices
\footnote{Submitted to J. Phys. A: Math. Gen.}}
\date{24th May 2001}
\author{
R. \'Alvarez-Nodarse${}^\dag{}^\ddag$ and R. S. Costas-Santos${}^\dag$ \\[5mm]
${}^\dag$ Departamento de An\'alisis Matem\'atico.\\
Universidad de Sevilla. Apdo. 1160, E-41080 Sevilla\\
${}^\ddag$ Instituto Carlos I de F\'{\i}sica Te\'orica y
Computacional, \\
Universidad de Granada, E-18071 Granada, Spain
}

\begin{document}
\maketitle
\begin{abstract}
We study the factorization of the hypergeometric-type difference
equation of Nikiforov and Uvarov on nonuniform lattices. An explicit
form of the raising and lowering operators is derived and some
relevant examples are given.
\end{abstract}
\section{Introduction}
In this paper we will deal with the so-called factorization
method (FM) of the hyper\-ge\-o\-me\-tric-type  difference equations
on nonuniform lattices. \ The FM was already used by
Darboux \cite{dar} and Schr\"odinger \cite{sch1,sch2}
to obtain the solutions of differential equations,
and also by Infeld and Hull \cite{inhu} for finding analytical
solutions of certain classes of second order differential
equations.
Later on,  Miller extended it to difference equations \cite{mil2}
and $q$-differences --in the Hahn sense-- \cite{mil3}. For more
recent works see e.g. \cite{ata1,atfrwo,ban2,smi3,spvizh} and references
therein.
\\
The classical FM was based on the existence of a so-called
raising and lowering  o\-pe\-ra\-tors for the corresponding equation
that allows to find the explicit solutions in a very easy way.
Going further, Atakishiyev and coauthors \cite{ata1,atwo,atfrwo}
have found the dynamical symmetry algebra related with the FM and
the differential or difference equations.
Of special interest was the paper by Smirnov \cite{smi1}
in which the equivalence of the FM and the Nikiforov et all theory
\cite{nisuuv} was shown, furthermore this paper pointed out
that the aforementioned equivalence remains valid also for the
nonuniform lattices that was shown later on in \cite{smi2,smi3}.
In particular, in \cite{smi3} a detailed study of the FM and
its equivalence with the Nikiforov et al. approach to difference
equations \cite{nisuuv} have been established. Also, in \cite{ban1},
a special nonuniform lattice was considered. In fact, in \cite{ban1}
the author constructed the FM for  the Askey--Wilson polynomials
using basically the difference equation  for the polynomials. In the
present paper we will continue the  research of the nonuniform lattice
case. In fact, following the idea by Bangerezako \cite{ban1} for
the Askey--Wilson polynomials and Lorente \cite{lor2} for the
classical continuous and discrete cases, we will obtain the FM for
the general polynomial solutions of the hypergeometric difference
equation on the general quadratic nonuniform lattice
$x(s)=c_1 q^s+c_2 q^{-s}+c_3$. We will use, as it is already
suggested in \cite{atwo,smi1}, not the polynomial solutions
but the corresponding normalized functions which is more
natural and useful. In such a way, the method proposed here is the
generalization of \cite{ban1} and \cite{lor2} to the aforementioned
nonuniform  lattice.\\
The structure of the paper is as follows. In Section 2 we
present some well-known results on orthogonal polynomials
on nonuniform lattices \cite{atrasu,nisuuv,niuv2}, in section
3 we introduce the normalized functions and obtain some
of their properties such as the lowering
and raising operator that allow us, in Section 4, to obtain
the factorization for the second order difference equation satisfied
by such functions. Finally, in Section 5, some relevant
examples are worked out.

\section{Some basic properties of the $q$-polynomials}

Here, we will summarize some of the properties of the {\it
q-}polynomials useful for the rest of the work. For further
information see e.g. \cite{nisuuv}.\\[0.3cm]
 We will deal here with the
{\it second order difference  equation of  the hypergeometric type}
\begin{equation}
\begin{array}{c}
\displaystyle
\sigma(s) \frac{\Delta}{\Delta x(s-\frac{1}{2})}
\left[\frac{\nabla y(s)}{\nabla x(s)}\right] + \tau(s)
\frac{\Delta y(s)}{\Delta x(s)} + \lambda y(s) =0, \\
\\
\sigma(s)=\tilde{\sigma}(x(s)) - \frac{1}{2}\tilde{\tau}(x(s)) \Delta
x(s-\frac{1}{2}), \quad \tau(s)=\tilde{\tau}(x(s)),
\end{array}
\label{eqdif}
\end{equation}
where $\nabla f(s)=f(s)-f(s-1)$ and $\Delta f(s)=f(s+1)-f(s)$
denote the backward and forward finite difference derivatives,
respectively, $\tilde{\sigma}(x(s))$ and $\tilde{\tau}(x(s))$
are polynomials in $x(s)$ of
degree at most  2 and 1, respectively, and $\lambda$ is a
constant. In the following, we will use the following notation for
the coefficients in the power expansions in $x(s)$ of $\tsig(s)$ and
$\ttau(s)$
\bq
\tsig(s)\equiv\tsig[x(s)]=\frac{\tsig''}2
x^2(s)+\tsig'(0)x(s)+\tsig(0),\quad \quad
\ttau(s)\equiv\ttau[x(s)]= \ttau' x(s)+ \ttau(0).
\label{sig-tau-pot}
\eq
An important property of the above
equation is that the {\em k-order difference derivative} of a
solution $y(s)$ of \refe{eqdif}, defined by
$$ y_{k}(s)_q =
\frac{\Delta}{\Delta x_{k-1}(s)}\frac{\Delta}{\Delta x_{k-2}(s)} \dots
\frac{\Delta}{\Delta x(s)}  y(s)  \equiv \Delta^{(k)}   y(s) , \,\,
$$
also satisfies a difference equation of the hypergeometric type
\bq
\sigma(s) \frac{\Delta  }{\Delta x_k(s-\frac{1}{2}) }\left[ \frac{\nabla
y_{k}(s)_q}{ \nabla x_k(s) }  \right] + \tau_k(s) \frac{\Delta
y_{k}(s)_q }{ \Delta x_k(s)} + \mu_k y_{k}(s)_q = 0,
\label{q-difder}
\eq
where $x_k(s)=x(s+\mbox{\footnotesize
$\frac{k}{2}$})$ and \cite[page 62, Eq. (3.1.29)]{nisuuv} \bq
\tau_k(s)=\frac{\sigma(s+k)-\sigma(s)+\tau(s+k)\Delta
x(s+k-\frac{1}{2})} {\Delta x_{k-1}(s)}, \quad\quad \label{tau_n}
\mu_k =\lambda + \sum_{m=0}^{k-1}\frac{\Delta \tau_m(s)}{\Delta
x_m(s)}. \label{mu_n}
\eq
It is important to notice that the above
difference equations have polynomial solutions of the hypergeometric
type iff $x(s)$ is a function of the form \cite{atrasu,niuv2}
\begin{equation}
x(s)=c_1(q) q^{s}+c_2(q)q^{-s}+c_3(q)=c_1(q)[q^{s}+q^{-s-\mu} ]+c_3(q),
\label{red-gen}
\end{equation}
where $c_1$, $c_2$, $c_3$ and $q^{\mu}=\frac{c_1}{c_2}$ are constants
which, in general, depend on $q$  \cite{nisuuv,niuv2}.
For the above lattice, a straightforward calculation shows that
$\tau_k(s)$ is a polynomial of first degree in $x_k(s)$ of the
form (see e.g. \cite{atrasu})
\bq
\begin{array}{c}
\tau_k(s)=\ttau_k'
x_k(s)+\ttau_k(0),\quad\quad \dst \ttau_k'=
[2k]_q\frac{\tsig''}2 +\alpha_q({2k})\ttau',\\[3mm]
\dst\ttau_k(0) =\frac{c_3 \widetilde{\sigma}''}{2}(2[k]_q-[2k]_q)+
\widetilde{\sigma}'(0)[k]_q+ c_3\tau'(\alpha_q(k)-\alpha_q({2k}))+
\ttau(0)\alpha_q(k),
\end{array}
\label{tau_k-taylor}
\eq
where the {\it q-numbers} $[k]_q$ and $\alpha_q(k)$ are defined by
\bq
[k]_q=\frac{q^{\frac{k}2}-q^{-\frac{k}2}}{q^{\frac{1}2}-q^{-\frac{1}2}},
\quad\quad \alpha_q(k)=\frac{q^{\frac{k}2}+q^{-\frac{k}2}}{2},
\label{q-num}
\eq
and $[n]_q!$ are the {\it q-}factorials $[n]_q! = [1]_q[2]_q\cdots[n]_q$.\

Both difference equations \refe{eqdif} and \refe{q-difder} can be
rewritten in the symmetric form $$ \frac{\Delta}{\Delta
x(\mbox{$s-\frac{1}{2}$})}\left[\sigma(s)\rho(s)\frac{\nabla
y(s)}{\nabla x (s)}\right]+\lambda_n \rho(s)y(s)=0, $$ and $$
\frac{\Delta}{\Delta x_k(s-\frac{1}{2})}\left[\sigma(s)\rho_k(s)\frac{\nabla
y_k(s)} {\nabla x_k(s)} \right] + \mu_k \rho_k(s)y_k(s) =0, $$ where
$\rho(s)$ and $\rho_k(s)$ are the weight functions satisfying the
Pearson-type difference equations \bq \ba{c} \displaystyle
\frac{\bigtriangleup}{\Delta x(s-\frac{1}{2})} \left[\sigma(s) \rho(s)
\right]= \tau(s) \rho(s)  \, , \quad\quad \dst
\frac{\bigtriangleup}{\Delta x_k(s-\frac{1}{2})}  \left[\sigma(s)
\rho_k(s)\right]=\tau_k(s) \rho_k(s),
 \ea
\label{pearson-q}
\eq
respectively.
In \cite{nisuuv} it is shown that the polynomial solutions of
\refe{q-difder} (and so the polynomial solutions of \refe{eqdif})
are determined by the $q$-analogue of the Rodrigues formula on the
nonuniform lattices
\bq
\frac{\Delta}{\Delta x_{k-1}(s)}\cdots
\frac{\Delta}{\Delta x(s)} P_n(x(s))_q \equiv \Delta^{(k)} P_n(x(s))_q =
\frac{A_{n,k} B_n}{\rho_k(s)} \nabla_k^{(n)}  \rho_n(s),
\label{roddif-q}
\eq
where
$$\ba{l}
 \dst   \nabla_k^{(n)} f(s)=
\frac{\nabla}{\nabla x_{k+1}(s)}\frac{\nabla}{\nabla  x_{k+2}(s)}
\cdots \frac{\nabla}{\nabla  x_n(s)} \,f(s) .
\ea
$$
\bq\hspace{.1cm}
A_{n,k} = \frac{[n]_q!}{[n-k]_q!}
\prod_{m=0}^{k-1} \left\{ \alpha_q(n+m-1)\widetilde{\tau}'+[n+m-1]_q
\frac{\widetilde{\sigma}''}{2}\right\}
\label{A_nm-q}
\eq
Thus \cite[ page 66, Eq. (3.2.19)]{nisuuv}
\begin{equation}
 P_n(x(s))_q=\frac{B_n}{\rho(s)} \nabla^{(n)}\rho _n(s),\quad
\nabla^{(n)}\equiv  \frac{\nabla}
{\nabla x_1(s)}\frac{\nabla}{\nabla x_2(s)}
\cdots \frac{\nabla}{\nabla  x_n(s)},
\label{rodeq}
\end{equation}
where $\rho_n(s)=\dst\rho(s+n) \prod_{k=1}^n \sigma(s+k)$ and
\bq
\lambda_n =
-[n]_q\left\{\alpha_q(n-1)\widetilde{\tau}'+[n-1]_q
\frac{\widetilde{\sigma}''}{2}\right\}.
\label{lambda-q}\label{lambda}
\eq
In this paper we will deal with  orthogonal {\it q-}polynomials and
functions. It can be proven \cite{nisuuv}, by using the difference equation of
hypergeometric-type \refe{eqdif}, that if the boundary condition
\begin{equation}
\begin{array}{c} \displaystyle
\sigma(s) \rho(s) x^{k}(s-\frac{1}{2}) \Big|_{s=a,b} = 0,
\quad   \forall k \ge  0,
\end{array}\label{con-ortho}
\end{equation}
holds, then the polynomials $P_n(s)_q$ are orthogonal, i.e.,
\begin{equation}
\begin{array}{c} \displaystyle
\sum_{s = a }^{b-1} P_n(x(s))_q P_m(x(s))_q\rho(s) \Delta
x(s-\frac{1}{2}) = \delta_{nm} d_n^2, \quad \!\!s=a,a+1,\dots,b-1,
\end{array}
\label{norm}
\end{equation}
where $\rho(s)$ is a solution of the Pearson-type equation
(\ref{pearson-q}). In the special case of the linear exponential
lattice $x(s)=q^s$ the above relation can be written in terms of
the Jackson $q$-integral (see e.g. \cite{gara,kost})
$\int_{z_1}^{z_2} f(t)d_q t$, defined by
$$
\int_{z_1}^{z_2} f(t)d_q t=
\int_{0}^{z_2} f(t)d_q t -\int_{0}^{z_1} f(t)d_q t,
$$
where
$$
\int_{0}^{z} f(t)d_q t=z(1-q)\sum_{k=0}^\infty f(zq^k) q^k,\quad
0<q<1,
$$
as follows:
\begin{equation}
\begin{array}{c} \displaystyle
\int_{q^a}^{q^b} P_n(t)_q P_m(t)_q \omega(t) d_q t
 = \delta_{nm} q^{1/2} d_n^2, \quad t=q^s,\quad \omega(t)\equiv
\omega(q^t)=\rho(t).
\end{array}
\label{norm-q-int}
\end{equation}

Notice that the above boundary condition
\refe{con-ortho} is valid  for $k=0$. Moreover, if we assume that
$a$ is finite, then \refe{con-ortho} is fulfilled at $s=a$ providing
that $\sigma(a)=0$ \cite[\S 3.3, page 70]{nisuuv}.
In the following we will assume that this condition holds.
The squared norm in \refe{norm} is given by
\cite[Chapter 3, Section {\bf 3.7.2}, pag. 104]{nisuuv}
$$
d_n^2 =(-1)^n A_{n,n} B_n^2 \sum_{s=a}^{b-n-1}\rho_n(s)
\Delta  x_n(s-\frac{1}{2}).
$$
There is also a so-called continuous orthogonality. In fact, if
there exist a contour $\Gamma$ such that
\bq
\int_{\Gamma}\Delta[\rho(z)\sigma(z) x^k(z-\frac{1}{2})] \, dz = 0,\quad\quad
\forall k\geq 0,
\label{condicion-orto-con}
\eq
then \cite{nisuuv}
$$
\int_{\Gamma}  P_n(x(z))_q P_m(x(z))_q \rho(z) \Delta x(z-\frac{1}{2}) \, dz = 0,
\quad\quad n\neq m.
$$
A simple consequence of the orthogonality is the following
three term recurrence relation:
\bq
\label{ttrr}
 x(s)P_n(x(s))_q=\alpha_n P_{n+1}(x(s))_q+\beta_n P_n(x(s))_q+
\gamma_n P_{n-1}(x(s))_q,
\end{equation}
where $\alpha_n$, $\beta_n$ and $\gamma_n$ are constants. If
$P_n(s)_q=a_n x^n(s)+b_n x^{n-1}(s)+ \cdots\,\,\,\, ,$ then using
(\ref{ttrr}) we find
\bq
\alpha_n= \frac{a_n}{a_{n+1}} ,\quad
\beta_n = \frac{b_n}{a_{n}}-\frac{b_{n+1}}{a_{n+1}},\quad
\gamma_n= \frac{a_{n-1}}{a_{n}}\frac{d_{n}^2}{d_{n-1}^2}.
\label{coef-RRTT-q}
\eq
To obtain the explicit values of
$\alpha_n$, $\beta_n$ we will use the following lemma,
--interesting in its own right-- that can be proven by induction:
\begin{lema}
\label{delta-k-x-n}
$$
\Delta^{(k)} x^n(s)= \frac{[n]_q!}{[n-k]_q!} x^{n-k}_k(s)+
c_3 \left(n\frac{[n-1]_q!}{[n-k-1]_q!}-(n-k)\frac{[n]_q!}{[n-k]_q!}\right)
 x^{n-k-1}_k(s)+\cdots.
$$
\end{lema}
In the case $k=n-1$, it becomes \bq \ba{l}\dst \Delta^{(n-1)}
x^n(s)= [n]_q! x_{n-1}(s)+ c_3[n-1]_q!\left(n-[n]_q\right).
\ea\label{delta-n-1-x-n}
\eq
Now, using the Rodrigues formula \refe{roddif-q} for $k=n-1$,
$$
\Delta^{(n-1)} P_n(x(s))_q =
\frac{A_{n,n-1} B_n}{\rho_{n-1}(s)} \nabla_{n-1}^{(n)}  \rho_n(s)=
\frac{A_{n,n-1} B_n}{\rho_{n-1}(s)} \frac{\nabla}{\nabla x_n(s)}\rho_n(s),
$$
as well as the identities $\rho_n(s)=\rho_{n-1}(s+1)\sigma(s+1)$,
$x_n(s)=x_{n-1}(s+\frac{1}{2})$ and the Pearson equation \refe{pearson-q}
for $\rho_{n-1}(s)$, we find
$$
\Delta^{(n-1)} P_n(x(s))_q  = A_{n,n-1} B_n \tau_{n-1}(s).
$$
Thus
$$
a_n=\frac{A_{n,n-1} B_n \ttau_{n-1}'}{[n]_q!}= B_n
\prod_{k=0}^{n-1} \left\{ \alpha_q(n+k-1) \tilde{\tau}' +
[n+k-1]_q \frac{\tilde{\sigma}''}{2} \right\} \,\,,
$$
and
$$
\frac{b_n}{a_n}=\frac{[n]_q\ttau_{n-1}(0)}{\ttau'_{n-1}}+c_3([n]_q-n).
$$
So
$$
\alpha_n=\frac{B_n}{B_{n+1}}\frac{\alpha_q({n-1})
\ttau'\!+\![n-1]_q\frac{\tsig''}2}
{(\alpha_q({2n\!-\!1})\ttau'\!+\![2n\!-\!1]_q\frac{\tsig''}2)
(\alpha_q({2n})\ttau'\!+\![2n]_q\frac{\tsig''}2)}=
-\frac{B_n}{B_{n+1}}\!\frac{\lambda_n}{[n]_q}\!
\frac{[2n]_q}{\lambda_{2n}}\!
\frac{[2n+1]_q}{\lambda_{2n+1}}
$$
and
$$ \beta_n=\frac{[n]_q\ttau_{n-1}(0)}{\ttau'_{n-1}}-
\frac{[n+1]_q\ttau_{n}(0)}{\ttau'_{n}}+c_3([n]_q+1-[n+1]_q).
$$
Using the Rodrigues formula the following difference-recurrent
relation follows \cite{alar1,nisuuv}

$$
 \sigma(s)\frac{\nabla P_n(x(s))_q}{\nabla
x(s)}=\frac{\lambda_n}{[n]_q\tau_n'}
\left[\tau_n(s)P_n(x(s))_q-\frac{B_n}{B_{n+1}}P_{n+1}(x(s))_q\right],
$$
where $\tau_n(s)$ is given by \refe{tau_k-taylor}, where the identity
$\ttau_n'=\dst -\frac{\lambda_{2n+1}}{[2n+1]_q}$ has been used.

Then, using the explicit expression for the coefficient $\alpha_n$, we find
\begin{equation}
\label{d3} \sigma(s)\frac{\nabla P_n(x(s))_q}{\nabla
x(s)}=\frac{\lambda_n}{[n]_q} \frac{\tau_n(s)}{\tau_n'}P_n(x(s))_q-
\frac{\alpha_n\lambda_{2n}}{[2n]_q}P_{n+1}(x(s))_q .
\end{equation}
This equation defines a raising operator in terms of the backward
difference in the sense that we can obtain the polynomial
$P_{n+1}$ of degree $n+1$ from the lower degree  polynomial
$P_n$.\\
From the above equation and using the identity
$\nabla=\Delta-\nabla \Delta$, the second order difference
equation and the three terms recurrence relation we find
\cite{alar1} lowering-type operator: \bq \label{d4}
\begin{array}{c}\dst [\sigma(s)+\tau(s)\Delta
x(\mbox{$s-\frac{1}{2}$})]\frac{\Delta P_n(x(s))_q}{\Delta x(s)}
=\frac{\gamma_n \lambda_{2n}}{[2n]_q}P_{n-1}(x(s))_q + \\[0.4cm] \dst
\left[\frac{\lambda_n}{[n]_q}\frac{\tau_n(s)}{\tau'_n}-\lambda_n
\Delta x (\mbox{$s-\frac{1}{2}$})-\frac{\lambda_{2n}}{[2n]_q}
(x(s)-\beta_n)\right]
P_n(x(s))_q.
\end{array}
\end{equation}
The most general polynomial solution of the {\it q-}hypergeometric
equation \refe{eqdif} corresponds to the case
\begin{equation}
\sigma(s) = A\prod_{i=1}^4 [s-s_i]_q = C q^{-2s}\prod_{i=1}^4(q^s-q^{s_i}),
\qquad A,C,\mbox{not vanishing constants}
\label{sigma-gen}
\end{equation}
and  has the form \cite{niuv2}
\begin{equation}
\begin{array}{l}\displaystyle
P_n(s)_q= \displaystyle   D_n  \,\, {}_{4}
\mbox{\large{$\phi$}}_3 \left(\begin{array}{c} q^{-n},
q^{ 2\mu+n-1+\mbox{\tiny $\sum_{i=1}^4$} s_i}, q^{s_1-s},q^{s_1+s+\mu} \\
q^{s_1+s_2+\mu},q^{s_1+s_3+\mu}, q^{s_1+s_4+\mu} \end{array}
;\, q\, ,\, q \right),
\end{array}
\label{rep-serie-q}
\end{equation}
where $D_n$ is a normalizing constant and the basic hypergeometric
series ${}_p  \mbox{\large{$\phi$}}_q$ are defined by \cite{kost}
$$ \mbox{\small $\displaystyle
\begin{array}{l} \displaystyle _{r} \mbox{\large{$\phi$}}_p
\left(\!\begin{array}{c} {a_1,\dots,a_{r}} \\
 {b_1,\dots,b_p} \end{array} \,;\, q\, ,\, z \!\right)=
 \displaystyle \sum _{k=0}^{\infty}\frac{ (a_1;q)_k \cdots  (a_{r};q)_k}
{(b_1;q)_k  \cdots  (b_p;q)_k}\frac{z^k}{(q;q)_k}
\left[ (-1)^k q^{\frac k2(k-1)}\right]^{p-r+1} ,
\end{array} $}
$$
and
\begin{equation}
(a;q)_k= \prod_{m=0}^{k-1}(1-aq^m),
\label{sim-poc-q}
\end{equation}
is the  $q$-analogue of the Pochhammer symbol. Instances of such
polynomials are the Askey--Wilson polynomials, the $q$-Racah
polynomials and big $q$-Jacobi polynomials among others
\cite{kost,niuv2}.


\section{The orthonormal functions on nonuniform lattices}

In this section we will introduce a set of orthonormal
functions which are orthogonal with respect to the unit
weight \cite{atwo,smi1}
\bq
\varphi_n(s)=\sqrt{\rho(s)/d_n^{2}}P_n(x(s))_q,
\label{nor-fun}
\eq
e.g. for the case of discrete orthogonality we have
$$
\sum_{s_i=a}^{b-1}\varphi_{n}(s_i)\varphi_{m}(s_i)\Delta x(s_i-\frac{1}{2})=
\delta_{nm}.
$$

Next, we will establish several important properties of such
functions which generalize, to the nonuniform lattices, the
ones given in \cite{lor2}. In the following we will use
the notation $\Theta(s)=\sigma(s)+\tau(s)\Delta x
(\mbox{$s-\frac{1}{2}$})$.

First of all, inserting \refe{nor-fun} into (\ref{eqdif}),
(\ref{ttrr}), (\ref{d3}), (\ref{d4})
we obtain that they satisfy the following difference equation:
\bq \label{nd1}\hspace{-.5cm}
\begin{array}{l}
\dst\sqrt{\Theta(s)\sigma(s+1)}\frac{1}{\Delta x(s)}\varphi_n(s+1)+
\dst\sqrt{\Theta(s-1))\sigma(s)}\frac{1}{\nabla
x(s)}\varphi_n(s-1)-
\\[0.4cm]
\dst\left(\frac{\Theta(s)}{\Delta x(s)}+\frac{\sigma(s)}{\nabla
x(s)}\right)\varphi_n(s)+\lambda_n \Delta x
(\mbox{$s-\frac{1}{2}$})\varphi_n(s)=0,
\end{array}
\end{equation}
the three term recurrence relation:
\bq \label{nd2}
\alpha_n \frac{d_{n+1}}{d_n}
\varphi_{n+1}(s)+\gamma_n
\frac{d_{n-1}}{d_n}\varphi_{n-1}(s)+(\beta_n-x(s))\varphi_n(s)=0,
\end{equation}
the raising-type formula:
\bq \label{lmas}\dst L^+(s,n)\varphi_n(s)
=\alpha_n\frac{\lambda_{2n}}{[2n]_q}\frac{d_{n+1}}{d_n}\varphi_{n+1}(s),
\end{equation}
and the lowering-type formula:
\bq
\label{lmenos}
L^-(s,n)\varphi_n(s)=\gamma_n\frac{\lambda_{2n}}{[2n]_q}
\frac{d_{n-1}}{d_n}\varphi_{n-1}(s),
\end{equation}
where  the raising-type operator $L^+(s,n)$ and the lowering-type
operator $L^-(s,n)$ are given by
\bq
\label{nd3}
\begin{array}{ll}
\dst L^+(s,n) \equiv & \dst
\left[\frac{\lambda_n}{[n]_q}\frac{\tau_n(s)}{\tau'_n}
-\frac{\sigma(s)}{\nabla x(s)}\right] I
{\dst
+\sqrt{\Theta(s-1)\sigma(s)}\frac{1}{\nabla x(s)}E^-},
\end{array}
\end{equation}
and
\bq \label{nd4}
\begin{array}{rl}
\dst L^-(s,n) \equiv & \dst
\left[-\frac{\lambda_n}{[n]_q}\frac{\tau_n(s)}{\tau'_n}+
\lambda_n \Delta x (\mbox{$s-\frac{1}{2}$})+
\frac{\lambda_{2n}}{[2n]_q}(x(s)-\beta_n) -
 \frac{\Theta(s)}{\Delta x(s)}\right] I\\[5mm]
& \dst
+\sqrt{\Theta(s)\sigma(s+1)}\dst\frac{1}{\Delta x(s)}E^+ ,
\end{array}
\end{equation}
respectively. In the above formulas $E^-f(s)=f(s-1)$, $E^+f(s)=f(s+1)$
and $I$ is the identity operator.

Notice that the last two formulas have a remarkable property of
giving all the solutions $\varphi_n(s)$. In fact, from
(\ref{nd4}) setting $n=0$ and taking into account that
$\varphi_{-1}(s)\equiv0$ we can obtain $\varphi_0(s)$. Then,
substituting the obtained function in (\ref{nd3}), we can find
all the functions $\varphi_1(s)$,\dots, $\varphi_n(s)$,\dots\,.
\begin{prop}
The raising and lowering operators (\ref{nd3}) and (\ref{nd4})
are mutually adjoint.
\end{prop}
\noindent{\bf Proof:} The proof is straightforward. In fact using the
boundary condition and after some calculations we obtain, in the case of
discrete orthogonality, the expression
$$\ba{l}\dst
\sum_{s_i=a}^{b-1}\varphi_{n+1}(s_i)\left[\frac{[2n]_q}{\lambda_{2n}}L^+(s_i,n)
\varphi_n(s_i)\right]\Delta x(s_i-\mbox{$\frac12$})\\[5mm]
\dst=\sum_{s_i=a}^{b-1}\left[\frac{[2n+2]_q}{\lambda_{2n+2}}
L^-(s_i,n+1)\varphi_{n+1}(s_i)\right] \varphi_n(s_i)
\Delta x(s_i-\mbox{$\frac12$})=\alpha_n\frac{d_{n+1}}{d_n}.
\ea
$$
The other cases can be done in an analogous way. \hfill\BBox
\begin{prop}
The operator corresponding to the eigenvalue $\lambda_n$
in (\ref{nd1}) is self adjoint.
\end{prop}

\noindent{\bf Proof:} Again we will prove the result in the case
of discrete orthogonality. Using the orthogonality conditions
$\sigma(a)\rho(a)=\sigma(b)\rho(b)=0$ (which is a consequence of
\refe{con-ortho}), we can write
$$
\sum_{s_i=a}^{b-1}\varphi_n(s_i)\sqrt{\Theta(s_i-1)
\sigma(s_i)}\frac{1}{\nabla
x(s_i)} \varphi_l(s_i-1)\Delta x(s_i-\mbox{$\frac12$})
$$
$$
=\sum_{s_i'=a-1}^{b-2}\varphi_n(s_i'+1)
\sqrt{\Theta(s_i')\sigma(s_i'+1)}\frac{1}{\nabla
x(s_i'+1)}\varphi_l(s_i')\Delta x(s'_i+\mbox{$\frac12$})
$$
$$
=\sum_{s_i=a}^{b-1}\varphi_n(s_i+1)
\sqrt{\Theta(s_i)\sigma(s_i+1)}\frac{1}{\nabla
x(s_i+1)}\varphi_l(s_i)\Delta x(s_i+\mbox{$\frac12$})+
$$
$$
\varphi_n(a)\sqrt{\Theta(a-1)\sigma(a)}\frac{1}{\nabla
x(a)}\varphi_l(a-1)\Delta x(a-\mbox{$\frac12$})-
$$
$$ \varphi_n(b)
\sqrt{\Theta(b-1)\sigma(b)}\frac{1}{\nabla
x(b)}\varphi_l(b-1)\Delta x(b-\mbox{$\frac12$}),
$$
where in the last two sums we first take the operations $\Delta$
and $\nabla$, and then substitute the corresponding value: e.g.
$\Delta x(a)=x(a+1)-x(a)$.

Now, we use the fact that
$\varphi_n(s)=\sqrt{\rho(s)/d_n^{2}}P_n(x(s))_q$,
as well as the boundary conditions
$\sigma(a)\rho(a)=\sigma(b)\rho(b)=0$, so
$$
\sqrt{\Theta(a-1)\sigma(a)}\,
\varphi_n(a) \varphi_l(a-1)=
\sqrt{\Theta(b-1)\sigma(b)}
 \varphi_n(b)\varphi_l(b-1)=0.
$$
The other terms can be transformed in a similar way.
All these yield the expression
$$ \sum_{s_i=a}^{b-1}
\varphi_l(s_i)\left\{\sqrt{\Theta(s_i)\sigma(s_i+1)}\frac{1}{\Delta
x(s_i)}\varphi_n(s_i+1)\Delta x(s_i+\mbox{$\frac12$})+\right.
$$
$$\left.
\qquad\sqrt{\Theta(s_i-1)\sigma(s_i)}\frac{1}{\nabla
x(s_i)}\varphi_n(s_i-1)\Delta x(s_i-\mbox{$\frac12$})\right\}=
$$
$$
=\sum_{s_i=a}^{b-1}
\varphi_n(s_i)\left\{\sqrt{\Theta(s_i)\sigma(s_i+1)}\frac{1}{\Delta
x(s_i)}\varphi_l(s_i+1)\Delta x(s_i+\mbox{$\frac12$})+\right.
$$
$$
\left. \qquad\sqrt{\Theta(s_i-1)\sigma(s_i)}\frac{1}{\nabla
x(s_i)}\varphi_l(s_i-1)\Delta x(s_i-\mbox{$\frac12$})\right\},
$$
from where the proposition easily follows. \hfill \BBox

\section{Factorization of difference equation of hypergeometric
type on the nonuniform lattice}

We will define from (\ref{nd1}) the following operator
$$
H(s,n)\equiv\sqrt{\Theta(s-1)\sigma(s)}\frac{1}{\nabla x(s)}E^- +
\sqrt{\Theta(s)\sigma(s+1)}\frac{1}{\Delta x(s)}E^+-
$$ $$
\left(\frac{\Theta(s)}{\Delta
x(s)}+\frac{\sigma(s)}{\nabla x(s)}-\lambda_n \Delta x
(\mbox{$s-\frac{1}{2}$})\right) I.
$$
Clearly, the orthonormal functions satisfy
$$
H(s,n)\varphi_n(s)=0.
$$
Let us rewrite the raising and
lowering operators in the following way
$$
\begin{array}{l} \dst  L^+(s,n)=u(s,n) I +
\sqrt{\Theta(s-1)\sigma(s)}\frac{1}{\nabla x(s)}E^-,
\\[0.6cm] \dst L^-(s,n)=v(s,n) I +
\sqrt{\Theta(s)\sigma(s+1)}\frac{1}{\Delta x(s)}E^+,
\end{array}
$$
where, as before, $\Theta(s)=
\sigma(s)+\tau(s) \Delta x (\mbox{$s-\frac{1}{2}$})$, and
$$
\begin{array}{l} \dst
u(s,n)=\frac{\lambda_n}{[n]_q}\frac{\tau_n(s)}{\tau'_n}-
\frac{\sigma(s)}{\nabla x(s)}, \\[0.5cm] \dst
v(s,n)=-\frac{\lambda_n}{[n]_q}\frac{\tau_n(s)}{\tau'_n}+\lambda_n
\Delta x (\mbox{$s-\frac{1}{2}$})+\frac{\lambda_{2n}}{[2n]_q}(x(s)-\beta_n)-
\frac{\Theta(s)}{\Delta x(s)}.
\end{array}$$
\begin{prop}\label{prop-appendix}
The functions $u(s,n)$ and $v(s,n)$ satisfy $u(s+1,n)=v(s,n+1)$, or,
equivalently  $u(s+1,n-1)=v(s,n)$.
\end{prop}
The proof of the above proposition is straightforward but
cumbersome. We will include it in  appendix A.
If we now calculate
$$
\ba{l} \dst L^-(s,n+1)L^+(s,n)=
v(s,n+1)u(s,n)+\Theta(s)\sigma(s+1)\left(\frac{1}{\Delta
x(s)}\right)^2+ \\[0.3cm] \dst + u(s+1,n)\left\{
\sqrt{\Theta(s-1)\sigma(s)}\frac{1}{\nabla x(s)}E^- +
 \sqrt{\Theta(s)
\sigma(s+1)}\frac{1}{\Delta x(s)}E^+
\right\},
\ea
$$
and substitute the values for $u(s,n)$, $v(s,n)$ and $H(s,n)$ we get
$$
\ba{l}
L^-(s,n+1)L^+(s,n)= h^\mp(n)I+u(s+1,n)H(s,n),
\ea
$$
where the function
$$
\ba{l}h^\mp(n)=
\dst\left.\left(\frac{\lambda_n}{[n]_q}
\frac{\tau_n(s+1)}{\tau_n'}-\frac{\sigma(s+1)}{\nabla x
(s+1)}\right)\left(\frac{\lambda_n}{[n]_q}\frac{\tau_n(s)}{\tau_n'}-\lambda_n
\Delta x (\mbox{$s-\frac{1}{2}$})\right)+ \right.\\[5mm]
\dst\quad\quad\quad\quad\left.\dst
\frac{\lambda_n}{[n]_q}\frac{\tau_n(s+1)}{\tau_n'}
\frac{\Theta(s)}{\Delta x
(s)}\right.,
\ea
$$
is independent of $s$. In fact, applying the last equality to the
orthonormal function $\varphi_n(s)$ and taking into account (\ref{lmas})
and  (\ref{lmenos}),
$$
h^\mp(n)=\frac{\lambda_{2n}}{[2n]_q}\frac{\lambda_{2n+2}}{[2n+2]_q}\alpha_n
\gamma_{n+1}.
$$
Similarly,
$$
L^+(s,n-1)L^-(s,n)= h^{\pm}(n)I+u(s,n-1)H(s,n),
$$
where
$$
\ba{rl}
h^{\pm}(n)=& \dst
\left.\left(-\frac{\lambda_n}{[n]_q}\frac{\tau_n(s-1)}{\tau_n'}
+\frac{\lambda_{2n}}{[2n]_q}(x(s-1)-\beta_n)+\lambda_n \Delta x
(\mbox{$s-\frac{3}{2}$})\right)\right.\times \\[5mm]
& \dst\left(-\frac{\lambda_n}{[n]_q}\frac{\tau_n(s)}
{\tau_n'}+\frac{\lambda_{2n}}{[2n]_q}(x(s)-\beta_n)+
\frac{\sigma(s)}{\nabla x(s)}\right)-\\[5mm]
& \dst\left(-\frac{\lambda_n}{[n]_q}
\frac{\tau_n(s)}{\tau_n'}+
\frac{\lambda_{2n}}{[2n]_q}(x(s)-\beta_n)\right)
\left(\frac{\Theta(s-1)}{\Delta
x(s-1)}\right),
\ea
$$
is independent of $s$. Furthermore, applying the last expression to
the functions $\varphi_n(s)$, and taking into account (\ref{lmas})
and (\ref{lmenos}), we obtain
$$
h^\pm(n)=\frac{\lambda_{2n-2}}{[2n-2]_q}
\frac{\lambda_{2n}}{[2n]_q}\alpha_{n-1}\gamma_n.
$$
{\bf Remark:} Notice that $h^\pm(n+1)=h^\mp(n)$.\\
All the above results lead us to our main theorem:

\begin{theo}
The operator $H(s,n)$, corresponding to the hypergeometric
difference equation for orthonormal functions $\varphi_n(s)$, admits
the following  factorization --usually called the Infeld-Hull-type
factorization--
\bq \label{nd5}
u(s+1,n)H(s,n)=L^-(s,n+1)L^+(s,n)-h^{\mp}(n)I,
\end{equation}
and \bq u(s,n)H(s,n+1)=\label{nd6}L^+(s,n)L^-(s,n+1)-h^\mp(n)I,
\end{equation}
respectively.
\end{theo}

\noindent {\bf Remark:} Substituting in the above formulas
the expression $x(s)=s$ we obtain the corresponding results
for the uniform lattice cases (Hahn, Kravchuk, Meixner and
Charlier), considered before by several authors, see e.g.
\cite{atwo,lor2,smi1} and by taking appropriate limits
(see e.g. \cite{kost,nisuuv}), we can recover the classical
continuous case (Jacobi, Laguerre and Hermite).

\section{Applications to some $q$-normalized ortogonal functions}
For the sake of completeness we will apply the above results to several
families of orthogonal $q$-po\-ly\-no\-mials and their corresponding
orthonormal $q$-functions that are of interest and appear in
several branches of mathematical physics. They are the
Al-Salam \& Carlitz polynomials I and II, the big $q$-Jacobi
polynomials, the dual $q$-Hahn polynomials, the continuous
$q$-Hermite and the celebrated $q$-Askey--Wilson polynomials.

The main data for these polynomials are taken from the
nice survey \cite{kost} except the case of dual $q$-Hahn
polynomials \cite{alsm}. Nevertheless, they can be
obtained also from the general formulas given in Section 2.

Finally, let us point out that similar factorization
formulas were obtained by other authors, e.g. Miller
in \cite{mil3} considered the polynomials on the
linear exponential lattice and Bangerezako studied the
Askey--Wilson case. Our main aim in this section is to show
how our general formulas lead, in a very easy way, to the
needed factorization formulas of several families for
normalized functions --not polynomials--.

\subsection{The Al-Salam \& Carlitz functions I and II}

The Al-Salam \& Carlitz polynomials I (and II) appear
in certain models of  $q$-harmonic oscillator , see e.g.
\cite{assu1,atsu2,atsu1,nag}.
They are polynomials on the exponential lattice
$x(s)=q^s\equiv x$, defined  \cite{kost} by
$$
U_n^{(a)}(x;q)=(-a)^n q^{\mbox{\tiny
$\left(\!\!\!\begin{array}{c} n \\[-0.05cm] 2
\end{array} \!\!\!\right)$}}\, _2
\phi_1\left(\begin{array}{c|c} q^{-n},x^{-1} \\[-0.3cm] &
q;\dst \frac{qx}{a} \\[-0.2cm]0 \end{array} \right),
$$
and constitute an orthogonal family with the orthogonality relation
\refe{norm-q-int}
$$
\int_a^1 U_n^{(a)}(x;q)U_m^{(a)}(x;q)\omega(x)d_q x=d_n^2
\delta_{nm},
$$
where
$$
\omega(x)=(qx,a^{-1}qx;q)_\infty,\quad\mbox{and}\quad
d_n^2=(-a)^n (1-q)(q;q)_n (q,a,a^{-1}q;q)_\infty
q^{\mbox{\tiny$\left(\!\!\!\begin{array}{c} n \\[-0.05cm] 2
\end{array} \!\!\!\right)$}}.
$$
As usual, $(a_1,\cdots,a_p;q)_n=(a_1;q)_n\cdots (a_p;q)_n$, and
$(a;q)_\infty=\prod_{k=0}^{\infty}(1-aq^k)$.\\
\noindent They satisfy a difference equation of
the form \refe{eqdif} where
$$\sigma(x)=(x-1)(x-a),\
\tau(x)=\widetilde{\tau}(x)=\tau' x+\tau(0),\ \mbox{being}\ \dst
\tau'=\frac{q^{1/2}}{1-q},\ \tau(0)=\dst q^{1/2}\frac{1+a}{q-1}.
$$
The eigenvalues $\lambda_n$ and the coefficients of the TTRR are
given by
$$
\dst \lambda_n=[n]_q\frac{q^{1-n/2}}{q-1}\quad\mbox{and}\quad
\alpha_n=1,\quad \beta_n=(1+a)q^n, \quad \gamma_n=aq^{n-1}(q^n-1),
$$
respectively. In this case we have
$$ \dst
\widetilde{\sigma}''=1,\quad \widetilde{\sigma}'(0)=-\dst
\frac{a+1}{2}, \quad \widetilde{\sigma}(0)=a, \quad \tau_n'=\dst
\frac{q^{\frac{1}{2}-n}}{1-q}, \quad \tau_n(0)=\dst
q^{\frac{1-n}{2}}\frac{a+1}{q-1}.
$$
The corresponding normalized functions \refe{nor-fun} are
$$
\varphi_n(x)=\sqrt{\frac{(qx,a^{-1}qx;q)_\infty(-a)^nq^{\mbox{\tiny
$\left(\!\!\!\begin{array}{c} n \\[-0.05cm] 2
\end{array} \!\!\!\right)$}}}{(1-q)(q;q)_n(q,a,q/a;q)_\infty}}\, _2
\varphi_1\left(\begin{array}{c|c} q^{-n},x^{-1} \\[-0.3cm] &
q;\dst \frac{qx}{a} \\[-0.2cm]0 \end{array} \right).
$$
Defining now the Hamiltonian for these functions $\varphi_n(x)$
$$ \dst H(x,n)\!\!= \!\dst  \!
\frac{\sqrt{a(x\!-\!1)(x\!-\!a)}}{x(1-q^{-1})}E^-
\!\!+\frac{\sqrt{a(q x\!-\!1)(q x\!-\!a)}} {x(q-1)}E^+\!+
\!\left(\!\frac{q^{1-n}}{1\!-\!q}x\!+\!\frac{q(a\!+\!1)}{q-1}\!-\!
\frac{[2]_q}{k_q}x^{-1}\!\right)I,  $$
 and using that $\dst
u(x,n)=\frac{aq}{1-q}x^{-1}$, $\dst
v(x,n)=u(qx,n-1)=\frac{a}{1-q}x^{-1}$, thus $$ L^+(x,n)=\dst
u(x,n)I+q\frac{\sqrt{a(x-1)(x-a)}}{x(q-1)}E^-,\quad \mbox{where}
\quad E^- f(x)=f(q^{-1}x), $$ and $$ \dst
L^-(x,n)=v(x,n)I+\frac{\sqrt{a(qx-1)(qx-a)}} {x(q-1)}E^+, \quad
\mbox{where} \quad E^+f(x)=f(qx), $$ we have $$ \dst
L^-(x,n+1)L^+(x,n)=\frac{aq^{1-n}(q^{n+1}-1)}
{(q-1)^2}I+v(x,n+1)H(x,n), $$ and $$ \dst L^+(x,n-1)L^-(x,n)=
\frac{aq^{2-n}(q^n-1)}{(q-1)^2}I+u(x,n-1)H(x,n),
$$
which give the factorization formulas for the Al-Salam \&
Carlitz functions I. If we now taking into account that
(see \cite[p. 115]{kost})
$$
V_n^{(a)}(x;q)=U_n^{(a)}(x;q^{-1}),
$$
then, the factorization for the Al-Salam \& Carlitz functions II
$$
\varphi_n(s)=q^{\mbox{\tiny
$\left(\!\!\!\begin{array}{c} s \\[-0.05cm] 2
\end{array} \!\!\!\right)$}}
\sqrt{\frac{a^{s+n}(aq;q)_\infty q^{\mbox{\tiny $\left(\!\!\!\!\!
\begin{array}{c} n\!\!+\!\!1 \\[-0.05cm] 2 \end{array} \!\!\!\!\!\right)$}}}
{(q,aq;q)_s(1-q)(q;q)_n}}\, _2\phi_0\left(\begin{array}{c|c}
q^{-n},x \\[-0.35cm] & q;\dst \frac{q^n}{a} \\[-0.2cm] - \end{array} \right),
$$
follows from the factorization for the Al-Salam \& Carlitz
functions I simply by changing $q$ to $q^{-1}$.

\subsection{The big $q$-Jacobi functions}

Now we will consider the most general family of $q$-polynomials
on the exponential lattice, the so-called  big $q$-Jacobi polynomials,
that appear in the representation theory of the quantum algebras
\cite{vikl}. They were introduced by Hahn in 1949 and are defined
\cite{kost} by
$$
P_n(x;a,b,c;q)=\frac{(aq;q)_n (cq;q)_n}{(abq^{n+1};q)_n} \,
_3\phi_2\left(\begin{array}{c|c} q^{-n},abq^{n+1},x \\[-0.2cm] &
q;q \\[-0.2cm] aq,cq \end{array} \right),\quad\quad
x(s)=q^s\equiv x.
$$
They constitute an orthogonal family, i.e.,
$$
\int_{cq}^{aq}\omega(x) P_n(x;a,b,c;q) P_n(x;a,b,c;q) d_q x=d_n^2
\delta_{nm}
$$
where
$$
\omega(x)=\dst \frac{(a^{-1}x;q)_\infty
(c^{-1}x;q)_\infty}{(x;q)_\infty (bc^{-1}x;q)_\infty},
$$
$$
d_n^2=
\frac{aq(1-q)(q,c/a,aq/c,abq^2;q)_\infty }{(aq,bq,cq,abq/c;q)_\infty}
\frac{(1-abq)(q,bq,abq/c;q)_n(-ac)^{-n}q^{\mbox{\tiny
$-\left(\!\!\!\begin{array}{c} n \\[-0.05cm] 2
\end{array} \!\!\!\right)$}}}{(abq,abq^{n+1},abq^{n+1})_n}.
$$

They satisfy the difference equation (\ref{eqdif}) with
$$
\sigma(x)=q^{-1}(x-aq)(x-cq) \ \mbox{ and } \ \tau(x)=\widetilde{\tau}(x)=\tau' x+\tau(0),
$$
where
$$
\tau'=\frac{1-abq^2}{(1-q)q^{1/2}} \quad \mbox{and} \quad
\tau(0)=\dst q^{1/2}\frac{a(bq-1)+c(aq-1)}{1-q},
$$
and
$$
\quad \lambda_n=-q^{-n/2}[n]_q\frac{1-abq^{n+1}}{1-q}.
$$
They satisfy a TTRR, whose coefficients are
$$
 \alpha_n=1, \quad  \beta_n=1-A_n-C_n, \quad \gamma_n=C_n A_{n-1},
$$
where
$$
\ba{l}
\dst A_n=\frac{(1-aq^{n+1})(1-cq^{n+1})(1-abq^{n+1})}
{(1-abq^{2n+1})(1-abq^{2n+2})},
 \\[0.5cm]
C_n=\dst -acq^{n+1}\frac{(1-q^n)(1-bq^n)(1-abc^{-1}q^n)}
{(1-abq^{2n})(1-abq^{2n+1})}.
\ea
$$
Also, we have
$$
\ba{l} \dst \widetilde{\sigma}''=\frac{1+abq^2}{q}, \quad
\widetilde{\sigma}'(0)=-\dst \frac{abq+acq+a+c}{2}, \quad
\widetilde{\sigma}(0)=acq,\\[0.4cm] \tau_n'=\dst
\frac{q^{-n}-abq^{n+2}}{q^{1/2}(1-q)}, \quad \tau_n(0)=\dst
q^{\frac{1-n}{2}}\frac{a(bq^{1+n}-1)+c(aq^{1+n}-1)}{1-q}. \ea
$$
The normalized big $q$-Jacobi functions are defined by
$$
\ba{c}
\dst \varphi_n(s)=\sqrt{\frac{(x/a,x/c;q)_\infty(aq,bq,abq/c;q)_
\infty (abq,aq,aq,cq,cq;q)_n(-ac)^n}{(x,bx/c,c/a,aq/c,abq^2;q)_\infty
(1-q)aq(1-abq)(q,bq,abq/c;q)_n}}\times \\[0.4cm]  \dst
\,_3\phi_2\left(\begin{array}{c|c} q^{-n},abq^{n+1},x \\[-0.2cm] &
q;q \\[-0.2cm] aq,cq
\end{array} \right).
\ea
$$
The corresponding Hamiltonian is
$$
\begin{array}{ll}
\!\! H(x,n)=& \!\!\!\! \dst
\frac{\sqrt{a(x\!-\!q)(x\!-\!a q)(x\!-\!c
q)(bx\!-\!cq)}}{x(q\!-\!1)}E^-\!\!+\! \dst
q\frac{\sqrt{a(x\!-\!1)(x\!-\!a)(x\!-\!c)(b x\!-\!c)}}
{x(q\!-\!1)}E^++
\\[0.5cm] & \dst
\left(\frac{1+abq^{2n+1}}{q^n(1-q)}x-\frac{q(a+ab+c+ac)}{1-q}+
\frac{acq(q+1)} {1-q}x^{-1}\right)I.
\end{array}
$$
Furthermore,
$$
u(x,n)=\dst \frac{abq^{n+1}}{1-q}x+D_n-\frac{acq^2}{q-1}x^{-1}
\quad \mbox{and } \ \dst
v(x,n)=\frac{abq^{n+1}}{1-q}x+D_{n-1}-\frac{acq}{q-1}x^{-1},
$$
where
$$
\dst D_n=\frac{ab(ab+ac+a+c)q^{2n+3}-a(b+c+ab+bc)q^{n+2}}{(1-abq^{2n+2})(1-q)},
$$
thus
$$
L^+(x,n)=\dst u(x,n)I+
\frac{\sqrt{a(x\!-\!q)(x\!-\!a q)(x\!-\!c q)(b x\!-\!c q)}}{x(q-1)}E^-,
\ \mbox{where } \  E^- f(x)=f(q^{-1}x),
$$
and
$$
L^-(x,n)=\dst v(x,n)I+q\frac{\sqrt{a(x\!-\!1)(x\!-\!a)(x\!-\!c)
(b x\!-\!c)}}{x(q-1)}E^+, \ \mbox{ where } \ E^+f(x)=f(qx),
$$
so
$$
L^-(x,n+1)L^+(x,n)=\delta_{n+1}\gamma_{n+1}I+v(x,n+1)H(x,n),
$$
$$
L^+(x,n-1)L^-(x,n)=\delta_n \gamma_n I+u(x,n-1)H(x,n),
$$
where
$$
\delta_n=\dst
\frac{(1-abq^{2n-1})(1-abq^{2n+1})}{q^{2n-1}(q-1)^2}.
$$
The above formulas are the factorization formulas
for the family of the big $q$-Jacobi normalized functions.\\
Since all discrete $q$-polynomials on the exponential lattice
$x(s)=c_1 q^s+c_3$ ---the so called, $q$-Hahn class--- can be obtained
from the big $q$-Jacobi polynomials by a certain limit process
(see e.g. \cite{alme,kost}, then from the above formulas we can obtain
the factorization formulas for the all other cases in the $q$-Hahn tableau.
Of special interest are the $q$-Hahn polynomials and the big $q$-Laguerre
polynomials, which are particular cases of the big $q$-Jacobi polynomials
when $c=q^{-N-1}$, $N=1,2,\dots$, and $c=0$, respectively.

\subsection{The $q$-dual-Hahn functions}
In this section we will deal with  the $q$-dual-Hahn polynomials,
introduced in \cite{alsm,niuv2} and closely related with the Clebsh-Gordon
coefficients of the $q$-algebras $SU_q(2)$ and $SU_q(1,1)$
\cite{alsm}. They are defined on the lattice
$x(s)=[s]_q[s+1]_q$ by
$$
W_n^c(x(s);a,b)_q=
\frac{(-1)^n (q^{a-b+1};q)_n (q^{a+c+1};q)_n}{q^{n/2(3a-b+c+1+n)}k_q^n
 (q;q)_n} \, _3\phi _2 \left(\begin{array}{c|c} q^{-n},q^{a-s},
q^{a+s+1} \\[-0.2cm] & q;q \\[-0.2cm] q^{a-b+1},q^{a+c+1}
\end{array} \right),
$$
and satisfy a discrete orthogonality \refe{norm} with
respect to the weight function
$$
\rho(s)=\dst
\frac{q^{\frac{1}{2}((b-1)^2-(2s-1)(a+c))}}
{(1-q)^{2(a+c-b)+1}}\frac{(q^{s-a+1},q^{s-c+1},q^{s+b+1},
q^{b-s};q)_\infty}{(q,q,q^{s+a+1},q^{s+c+1};q)_\infty},
$$
where $-\frac{1}{2}\le a<b-1, \quad |c|<a+1,$
and for this weight function the norm is
$$
d_n^2=\dst \frac{q^{\frac{1}{4}(-4ab-4bc+6a+6c-8b+6+4n(a+c-2b)-
n^2+17n+2b^2)}}{(1-q)^{2(a+c-b+1)+3n}}\dst
 \frac{(q^{b-c-n},q^{b-a-n};q)_\infty}{[n]_q!(q,q^{a+c+n+1};q)_\infty}.
$$
These polynomials satisfy a TTRR (\ref{ttrr}) with
$$
\ba{l} \alpha_n\!=\!1, \\[0.3cm]
\beta_n\!\!=\!\!q^{\frac{1}{2}(2n\!-\!b\!+\!c\!+\!1)}[b\!-\!a\!-\!n\!+\!1]_q[a\!+\!c\!+\!n\!+\!1]_q\!+\!
q^{\frac{1}{2}(2n\!+\!2a\!+\!c\!-\!b\!+\!1)}[n]_q[b\!-\!c\!-\!n]_q\!+\![a]_q[a\!+\!1]_q,\\[0.3cm]
\gamma_n\!=\!q^{2n+c+a-b}[a+c+n]_q[b-a-n]_q [b-c-n]_q[n]_q, \ea
$$
and the second order difference equation~(\ref{eqdif}),
whose eigenvalues are $\lambda_n=[n]_q q^{\frac{1}{2}-\frac{n}{2}}$ and
$$
\sigma(s)=q^{\frac{1}{2}(s+c+a-b+2)}[s-a]_q[s+b]_q[s-c]_q
\quad \mbox{and} \quad \tau(x)=\widetilde{\tau}(x)=\tau' x+\tau(0),
$$
with $\tau'=-1$ and
$\tau(0)=q^{\frac{1}{2}(a-b+c+1)}[a+1]_q[b-c-1]_q+q^{\frac{1}{2}(c-b+1)}[b]_q[c]_q.$
\\[0.3cm]
Also we will need the values
$$
\ba{l}
\widetilde{\sigma}''\ =\ k_q,\quad\dst\widetilde{\sigma}'(0)=\frac{1}{2k_q}
(2[2]_q-q^{\frac{1}{2}-b}-q^{\frac{1}{2}+a}-q^{\frac{3}{2}+a+c-b}-q^{\frac{1}{2}+c}),\\[0.3cm]
\widetilde{\sigma}(0)\!=\!\dst \frac{1}{2k_q^3}(2q^{1\!+\!a\!-\!b}\!+\!q^{-1}+\!q\!+
\!2q^{1\!+\!c\!-b}\!+2q^{1\!+\!a\!+\!c}\!-
\!(1\!+\!q)(q^{-b}\!+\!q^a\!+\!q^c\!+\!q^{1\!+\!a\!+\!c\!-\!b})),
\\[0.4cm] \tau_n'\!=\!-q^{-n}, \quad
\tau_n(0)\!=\!q^{\frac{1}{2}(c-b-n+1)}
[c\!+\!\frac{n}{2}]_q [b\!-\!\frac{n}{2}]_q\!+\!q^{\frac{1}{2}(a\!+\!c\!-\!b\!+\!1
\!-\!\frac{n}{2})}[a\!+\!\frac{n}{2}\!+\!1]_q[b\!-\!c\!-\!n\!-\!1]_q,
\ea
$$
In this case, the Hamiltonian, associated with the $q$-dual Hahn
normalized functions
$\sqrt{\rho(s)/d^{2}_n} W_n^c(x(s);a,b)_q$, is
$$
\ba{l} H(s,n)=\dst q^{\frac{1}{2}(c+a-b+2)}\frac{\sqrt{([s+1]_q^2
-[a]_q^2)([b]^2_q-[s+1]_q^2)([s+1]_q^2-[c]_q^2)}}{[2s+2]_q}E^+ +\\[0.4cm]
\dst q^{\frac{1}{2}(c+a-b+2)}\frac{\sqrt{([s]_q^2-[a]_q^2)([b]^2_q-[s]_q^2)
([s]_q^2-[c]_q^2)}}{[2s]_q}E^- -q^{\frac{1}{2}-\frac{n}{2}}[n]_q
[2s+1]_qI+\\[0.4cm]
\dst  q^{\frac{1}{2}(c+a-b+2)}\left(\frac{[s-a]_q[s+b]_q+[s-c]_q}{[2s]_q}-
\frac{[s+1-a]_q[s+1+b]_q[s+1-c]_q}{[2s+2]_q}\right)I,\ea
$$
where $E^+ f(s)=f(s+1)$ and  $E^- f(s)=f(s-1)$.
Then, using that
$$\ba{l} u(s,n)=q^{\frac{1}{2}-\frac{n}{2}}x(s+n/2)-
q^{\frac{1}{2}+\frac{n}{2}}(q^{\frac{1}{2}(c-b-n+1)}
[c+\frac{n}{2}]_q [b-\frac{n}{2}]_q+\\[0.4cm]
q^{\frac{1}{2}(a+c-b+1-\frac{n}{2})}[a+\frac{n}{2}+1]_q[b-c-n-1]_q)-q^{\frac{1}{2} (s+c+a-b+2)}
\dst \frac{[s-a]_q[s+b]_q[s-c]_q}{[2s]_q},\ea$$
and taking into account that
$v(s,n)=u(s+1,n-1)$, we find
$$
L^+(s,n)=u(s,n)I+q^{\frac{1}{2}(c+a-b+2)}
\frac{\sqrt{([s]_q^2-[a]_q^2)([b]^2_q-[s]_q^2)
([s]_q^2-[c]_q^2)}}{[2s]_q}E^-,
$$
and
$$
L^-(s,n)=v(s,n)I+q^{\frac{1}{2}(c+a-b+2)}
\frac{\sqrt{([s+1]_q^2-[a]_q^2)([b]^2_q-[s+1]_q^2)
([s+1]_q^2-[c]_q^2)}}{[2s+2]_q}E^+.
$$
Thus
$$
L^-(s,n+1)L^+(s,n)=q^{-2n}\gamma_{n+1}I+v(s,n+1)H(s,n),
$$
and
$$
L^+(s,n-1)L^-(s,n)=q^{-2n+2}\gamma_n I+u(s,n-1)H(s,n),
$$
are the factorization formulas for the $q$-dual Hahn normalized functions.

\subsection{The Askey--Wilson functions}
Finally we will consider the family of Askey--Wilson polynomials.
They are polynomials on the lattice
$x(s)=\frac{1}{2}(q^s+q^{-s})\equiv x$, defined by \cite{kost}
$$
p_n(x(s);a,b,c,d)=\dst \frac{(ab;q)_n (ac;q)_n (ad;q)_n}{a^n}{}_4 \phi_3
\left( \begin{array}{c|c} q^{-n},
q^{n-1}abcd,ae^{-i\theta},ae^{i\theta} \\[-0.2cm] & q;q
\\[-0.2cm] ab,ac,ad \end{array} \right),
$$
i.e., they correspond to the general case \refe{rep-serie-q} when
$q^{s_1}=a$, $q^{s_2}=b$, $q^{s_3}=c$, $q^{s_4}=d$.
Their orthogonality relation is of the form
$$
\int_{-1}^1 \omega(x) p_n(x;a,b,c,d) p_m(x;a,b,c,d)
\sqrt{1-x^2}\kappa_q dx =
\delta_{nm} d_n^2,\quad\quad q^s=e^{i\theta},\quad x=\cos\theta,
$$
where
$$
\omega(x)=
\frac{h(x,1)h(x,-1)h(x,q^{\frac12})h(x,-q^{\frac12})}
{\dst2\pi\kappa_q(1-x^2)h(x,a)h(x,b)h(x,c)h(x,d)}, \quad
h(x,\alpha)=\prod_{k=0}^\infty [1-2\alpha xq^k+\alpha^2q^{2k}],
$$
and the norm is given by
$$
d_n^2=\frac{(abcdq^{n-1};q)_n(abcdq^{2n};q)_\infty}
{(q^{n+1},abq^n,acq^n,adq^n, bcq^n,bdq^n,cdq^n;q)_\infty}.
$$
The Askey--Wilson polynomials satisfy the difference equation
(\ref{eqdif}) with
$$
\sigma(s)=-q^{-2s+1/2}\kappa_q^2
(q^s-a)(q^s-b)(q^s-c)(q^s-d),\qquad
\kappa_q=(q^{\frac{1}{2}}-q^{-\frac{1}{2}})
$$
and $\tau(x)=\widetilde{\tau}(x)=\tau' x+\tau(0)$,
where
$$
\tau'=4(q-1)(1-abcd),\ \tau(0)=2(1-q)(a+b+c+d-abc-abd-acd-bcd).
$$
Furthermore, they satisfy the TTRR \refe{ttrr} with
coefficients
$$
\alpha_n=1, \quad \beta_n=\frac{a+a^{-1}-(A_n+C_n)}{2},
\quad \gamma_n=\frac{C_nA_{n-1}}4,
$$
where $A_n$, $C_n$ are defined by
$$
\ba{l} \dst
A_n=\frac{(1-abq^n)(1-acq^n)(1-adq^n)(1-abcdq^{n-1})}
{a(1-abcdq^{2n-1})(1-abcdq^{2n})},\\[0.3cm]
C_n=\dst \frac{a(1-q^n)(1-bcq^{n-1})(1-bdq^{n-1})(1-cdq^{n-1})}
{(1-abcdq^{2n-2})(1-abcdq^{2n-1})},
\ea
$$
and whose eigenvalues are
$\lambda_n= 4  q^{-n+1} (1-q^n)(1-abcdq^{n-1})$.
In addition,  we have
$$
\ba{l} \dst
\widetilde{\sigma}''=-4(q-1)^2(1+abcd)q^{-1/2}, \\[0.3cm] \dst
\widetilde{\sigma}'(0)=(q-1)^2(a+b+c+d+abc+abd+acd+bcd)q^{-1/2},\\[0.3cm]
\dst
\widetilde{\sigma}(0)=(q-1)^2(1-ab-ac-ad-bc-bd-cd+abcd)q^{-1/2},
\\[0.3cm] \dst \tau_n'=4q^{-n}(q-1)(1-abcdq^{2n}),\\[0.3cm] \dst
\tau_n(0)=2(q-1)(-a-b-c-d+(abc+abd+acd+bcd)q^n)q^{-n/2}. \ea
$$
Defining now the normalized functions (see \refe{norm-q-int})
$\sqrt{\omega(x)/d^{2}_n} p_n(x;a,b,c,d)$, the corresponding
Hamiltonian $H(s,n)$ is
$$\ba{rl}
H(s,n)=& \dst \frac{2q^{3/2}}{[2s-1]_q} G(s,a,b,c,d)
E^- +
\frac{2q^{3/2}}{[2s+1]_q} G(s+1,a,b,c,d)
E^+ \, +  \\[5mm]
& \dst 2\Bigg(q^{-2s+1/2}\frac{\prod_{i=1}^4(1-q^{s_i+s})}
{[2s+1]_q}+q^{-2s+1/2}\frac{\prod_{i=1}^4
(q^s-q^{s_i})}{[2s-1]_q}+\\[5mm]
\dst& \qquad q^{-n+1}\kappa_q^2
(1-q^n)(1-abcdq^{n-1})[2s]_q\Bigg)I
\ea
$$
where
$$
G(s,a,b,c,d)=\sqrt{\prod_{i=1}^4
(1-2q^{s_i} q^{-1/2} x(s-{1}/{2})+q^{-1}q^{2s_i})},\quad
$$
We now define
$$
u(s,n)=D_n x_n(s)+D_nE_n+\dst q^{-2s+1/2}
\frac{(q^s-a){(q^s-b)}{(q^s-c)}{(q^s-d)}}{[2s-1]_q}
$$
where
$$
D_n=-4q^{-n/2+1/2}(q-1)(1-abcdq^{n-1}). \quad
$$
$$
E_n=\dst
\frac{(-a-b-c-d+(abc+abd+acd+bcd)q^n)q^{n/2}}{2(1-abcdq^{2n})}.
$$
Taking into account that $v(s,n)=u(s+1,n-1)$, we find
$$
L^+(s,n)=u(s,n)I+ \frac{2q^{3/2}}{[2s-1]_q}
G(s,a,b,c,d) E^-,
$$
$$
L^-(s,n)=v(s,n)I+ \frac{2q^{3/2}}{[2s+1]_q} G(s+1,a,b,c,d)
E^+,
$$
where $E^-f(s)=f(s-1)$  and $E^+f(s)=f(s+1)$.
Thus,
$$ L^-(s,n+1)L^+(s,n)=D_{2n}D_{2n+2}\gamma_{n+1} I+v(s,n+1)H(s,n),
$$
and
$$
L^+(s,n-1) L^-(s,n)=D_{2n-2}D_{2n}\gamma_{n}
I+u(s,n-1)H(s,n),
$$
which is the factorization formula for the Askey--Wilson
functions. \\
To conclude this paper let us consider the special case of
Askey--Wilson polynomials when $a=b=c=d=0$, i.e., the
continuous $q$-Hermite polynomials
$$
H_n(x|q)=2^{-n} e^{i n\theta} {}_2\phi_0
\left(\begin{array}{c|c} q^{-n},
0 \\[-0.2cm] & q;q^n e^{-2i\theta}
\\[-0.2cm] \mbox{---} \end{array} \right),\qquad x=\cos\theta.
$$
These polynomials are closely related with the $q$-harmonic oscilator
model introduced by Biedenharn \cite{bie} and Macfarlane
\cite{mac}, as it was pointed out in \cite{atsu1},
where the factorization for the continuous $q$-Hermite
polynomials were considered first. If we substitute $a=b=c=d=0$
in the above formulas, we obtain the factorization
for the $q$-Hermite functions
$$
\varphi_n(x)= \sqrt{
\frac{h(x,1)h(x,-1)h(x,q^{1/2})h(x,-q^{1/2})
(q^{n+1};q)_\infty}
{2\pi\kappa_q(1-x^2)}}\,H_n(x|q).
$$
In fact, since for continuous $q$-Hermite polynomials
$$
\sigma(s)=-\kappa_q^2 q^{2s+1/2},\quad
\tau(s)=4(q-1)x(s),\quad\lambda_n=4q^{-n+1}(1-q^n),
$$
and the coefficients for the three-term recurrence relation
are $\alpha_n=1$, $\beta_n=0$, $\gamma_n=(1-q^n)/4$, then
we obtain
$$
H(s,n)=\frac{2q^{3/2}}{[2s-1]_q}E^- +
\frac{2q^{3/2}}{[2s+1]_q}E^+ +2\!
\left(\!\frac{q^{-2s+1/2}}{[2s+1]_q}\!+\!
\frac{q^{2s+1/2}}{[2s-1]_q}\!-\!q^{-n+1}\kappa_q^2(1\!-\!q^n)
[2s]_q\! \right)\!I,
$$
$$
L^+(s,n)=\left(-4q^{-n/2+1/2}(q-1)x(s+n/2)+
\frac{q^{2s+1/2}}{[2s-1]_q}\right)I +
\frac{2q^{3/2}}{[2s-1]_q}E^-,
$$
$$
L^-(s,n)=\left(-4q^{-n/2+1}(q-1)x(s+n/2+1/2)+
\frac{q^{2s+5/2}}{[2s+1]_q}\right)I +
\frac{2q^{3/2}}{[2s+1]_q}E^-
$$
and $h^{\pm}(n)=4\kappa_q^2 q^{-2n+1}(1-q^n)$.

\subsection*{Acknowledgements}
The authors thank N. Atakishiyev and Yu. F. Smirnov for interesting
discussions and remarks that allowed us to improve this paper substantially,
as well as the referees for their remarks. The work has been partially
supported by the Ministerio de Ciencias y Tecnolog\'{\i}a of Spain under
the  grant BFM-2000-0206-C04-02, the Junta de Andaluc\'{\i}a
under grant FQM-262 and the European proyect INTAS-2000-272.

\bigskip\bigskip

\noindent\appendix{\bf\large Appendix A}

{\small

\vskip .5cm
\noindent Here, for the sake of completeness, we will prove
Proposition \ref{prop-appendix}, by showing that $u(s+1,n)-v(s,n+1)=0$.
To do that, we start with computing the difference

$$\begin{array}{l} \dst u(s+1,n)-v(s,n+1)=
 \frac{\lambda_n}{[n]_q}
\frac{\tau_n(s+1)}{\tau_n'}- \frac{\Delta \sigma(s)}{\Delta x (s)}+
 \\[0.4cm] \dst \frac{\lambda_{n+1}}{[n+1]_q}\frac{\tau_{n+1}(s)}
{\tau'_{n+1}}-\lambda_{n+1} \Delta x \mbox{$\left( s-\frac{1}{2}
\right)$}-\frac{\lambda_{2n+2}}{[2n+2]_q}(x(s)-\beta_{n+1})+
\frac{\tau(s)\Delta x \mbox{$\left( s-\frac{1}{2}
\right)$}}{\Delta x(s)}. \ea $$ Now we use the expansion
$\tau_n(s+1)=\tau_n' x_n(s+1)+\tau_n(0)$. Since $$ \ba{l} \dst
\frac{\Delta (x^2(s))}{\Delta
x(s)}=\frac{x^2(s+1)-x^2(s)}{x(s+1)-x(s)}=x(s+1)+x(s)=C_1
q^s(q+1)+C_2q^{-s}(q^{-1}+1)+2C_3= \\[0.4cm]
(C_1q^{s+\frac{1}{2}}+C_2q^{-s-\frac{1}{2}})[2]_q+2C_3=[2]_qx_1(s)+(2-[2]_q)C_3,
\ea $$ $$\ba{l}\dst x(s)\Delta
x(s-\frac{1}{2})=x(s)(C_1q^{s-\frac{1}{2}}(q-1)+C_2q^{-s+\frac{1}{2}}(q^{-1}-1))=
x(s)(C_1q^s-C_2q^{-s})k_q= \\[0.4cm]
\dst(C_1^2q^{2s}-C_2^2q^{-2s})k_q+ C_3(C_1q^s-C_2q^{-s})k_q, \ea
$$ where $k_q=q^{\frac{1}{2}}-q^{-\frac{1}{2}}$, $$ \dst \frac{\Delta}{\Delta
x(s)}\Big(x(s)\Delta x (\mbox{$s-\frac{1}{2}$})\Big)=\left(\frac{(C_1^2
q^{2s+1}+C_2^2
q^{-2s-1})[2]_q+C_3(C_1q^{s+\frac{1}{2}}+C_2q^{-s-\frac{1}{2}})}
{C_1q^{s+\frac{1}{2}}-C_2q^{-s-\frac{1}{2}}}\right)k_q, $$ and $$ \dst
\frac{\Delta}{\Delta x(s)}\Big( \Delta
x(s-\frac{1}{2})\Big)=\frac{\Delta}{\Delta x(s)}\Big((C_1
q^{s}-C_2q^{-s})k_q\Big)=\frac{C_1q^{s+\frac{1}{2}}+C_2q^{-s-\frac{1}{2}}}
{C_1q^{s+\frac{1}{2}}-C_2q^{-s-\frac{1}{2}}}k_q. $$ Then $$ \ba{l} \dst
\frac{\Delta \sigma(s)}{\Delta x(s)}=\frac{\Delta}{\Delta x(s)}
\left(\widetilde{\sigma}(s)-\frac{1}{2}\widetilde{\tau}(s)\Delta
x(\mbox{$s-\frac{1}{2}$})\right)=\\[0.4cm] \dst
=\frac{\Delta}{\Delta x(s)} \left(
\frac{\widetilde{\sigma}''}{2}x^2(s)+\widetilde{\sigma}'(0)x(s)+\widetilde{\sigma}(0)-
\frac{1}{2}\left(\tau' x(s)+\tau(0)\right)\Delta x
(\mbox{$s-\frac{1}{2}$}) \right)=\\[0.4cm] \dst
\frac{\widetilde{\sigma}''}{2}\left([2]_qx_1(s)+(2-[2]_q)C_3\right)
+\widetilde{\sigma}'(0)-\frac{1}{2}\tau(0)\left(\frac{C_1q^{s+\frac{1}{2}}+C_2q^{-s-\frac{1}{2}}}
{C_1q^{s+\frac{1}{2}}-C_2q^{-s-\frac{1}{2}}}\right)k_q-\\[0.4cm] \dst
 \frac{1}{2}\tau'
\left(\frac{[2]_q(C_1^2q^{2s+1}+C_2^2
q^{-2s-1})+C_3(C_1q^{s+\frac{1}{2}}+C_2q^{-s-\frac{1}{2}})}
{C_1q^{s+\frac{1}{2}}-C_2q^{-s-\frac{1}{2}}}\right) k_q.
 \ea $$
This yields for $u(s+1,n)-v(s,n+1)$ the expression
$$
\ba{l} \dst =\left[\frac{\lambda_n}{[n]_q}
x_n(s+1)+\frac{\lambda_n}{[n]_q}\frac{\tau_n(0)}{\tau_n'}\right]
-\left[\frac{\widetilde{\sigma}''}{2}[2]_qx_1(s)+\frac{C_3}{2}(2-[2]_q)
\widetilde{\sigma}''+ \widetilde{\sigma}'(0)-\right.
\ea
$$
$$\ba{l}
\dst
\frac{\widetilde{\tau}'}{2} \left(\frac{[2]_q(C_1^2q^{2s+1}+C_2^2
q^{-2s-1})}{C_1q^{s+\frac{1}{2}}-C_2q^{-s-\frac{1}{2}}}+ \frac{C_3
x_1(s)-C_3^2} {C_1q^{s+\frac{1}{2}}-C_2q^{-s-\frac{1}{2}}} \right)k_q-
\\[0.4cm] \dst \left. \frac{\tau(0)}{2}\left(
\frac{x_1(s)-C_3}{C_1q^{s+\frac{1}{2}}-C_2q^{-s-\frac{1}{2}}}\right)k_q
\right]+
\frac{\lambda_{n+1}}{[n+1]_q}\frac{\tau_{n+1}(s)}{\tau'_{n+1}}
 -\lambda_{n+1} \Delta x \mbox{$\left( s-\frac{1}{2}
\right)$}-\\[0.6cm] \dst
\frac{\lambda_{2n+2}}{[2n+2]_q}\left[C_1q^s+C_2q^{-s}+C_3-\frac{[n+1]_q\tau_n(0)}{\tau'_n}
\right.+\\[0.6cm] \dst \left.
\frac{[n+2]_q\tau_{n+1}(0)}{\tau_{n+1}'}-C_3(1+[n+1]_q-[n+2]_q)\right]+
\frac{\tau(s)\Delta x \mbox{$\left( s-\frac{1}{2}
\right)$}}{\Delta x(s)}. \ea $$
 Next, we expand $\Delta x_n(s)$ and
$\frac{\widetilde{\sigma}''}{2}[2]x_1(s)$, make some straightforward
calculations and  use the identities: $$
\frac{\lambda_n}{[n]_q}\frac{\tau_n(0)}{\tau_n'}+[n+1]_q\frac{\lambda_{2n+2}}{[2n+2]_q}
\frac{\tau_n(0)}{\tau_n'}=\left(\frac{\lambda_n}{[n]_q}+[n+1]_q\frac{\lambda_{2n+2}}{[2n+2]_q}\right)
\frac{\tau_n(0)}{\tau'_n}=-[n+2]_q\tau_n(0), $$ $$
\frac{\lambda_{n+1}}{[n+1]_q}\frac{\tau_{n+1}(s)}{\tau_{n+1}'}-[n+2]_q
\frac{\lambda_{2n+2}}{[2n+2]_q}\frac{\tau_{n+1}(0)}{\tau_{n+1}'}=[n+1]_q\tau_{n+1}(0)+
\frac{\lambda_{n+1}}{[n+1]_q}x_{n+1}(s), $$ as well as $$\ba{l} \dst
\frac{\lambda_n}{[n]_q}(C_1q^{s+1+\frac{n}{2}}+C_2q^{-s-1-\frac{n}{2}})-
\frac{\widetilde{\sigma}''}{2}[2]_q(C_1q^{s+\frac{1}{2}}+C_2q^{-s-\frac{1}{2}})-\frac{\lambda_{2n+2}}{[2n+2]_q}
(C_1q^s+C_2q^{-s})-\\[0.4cm]
\dst\lambda_{n+1}(C_1q^s-C_2q^{-s})k_q+\frac{1}{2}
\tau'(C_1q^{s+\frac{1}{2}}+C_2q^{-s-\frac{1}{2}})(q+q^{-1})=\frac{C_1q^s
\tau'}{2}(q^{n+\frac{1}{2}}+q^{\frac{1}{2}})+\\[0.4cm] \dst
\frac{C_1q^s\widetilde{\sigma}''}{2(q^{\frac{1}{2}}-q^{-\frac{1}{2}})}
(q^{n+\frac{1}{2}}-q^{\frac{1}{2}})+\frac{C_2q^{-s}\tau'}{2}(q^{-n-\frac{1}{2}}+q^{-\frac{1}{2}})+
\frac{C_2q^{-s}\widetilde{\sigma}''}{2(q^{\frac{1}{2}}-q^{-\frac{1}{2}})}(-q^{-n-\frac{1}{2}}+q^{-\frac{1}{2}})=
\\[0.4cm]
\dst =
-\frac{\lambda_{n+1}}{[n+1]_q}(C_1q^{s+\frac{n+1}{2}}+C_2q^{-s-\frac{n+1}{2}}),
\ea $$ we find $$\ba{l}
\dst=-\frac{\lambda_{n+1}}{[n+1]_q}\big(C_1q^{s+\frac{n+1}{2}}+C_2q^{-s-
\frac{n+1}{2}}\big)+C_3\frac{\lambda_n}{[n]_q}-[n+2]_q\tau_n(0)-
C_3\widetilde{\sigma}''-\widetilde{\sigma}'(0)+\frac{1}{2}
\tau'C_3k_q+\\[0.4cm]
\!\!\!\!\!\!\quad\dst\frac{1}{2}\tau(0)k_q\!+ \!
[n\!+\!1]_q\tau_{n+1}(0)\!+\!
\frac{\lambda_{n+1}}{[n\!+\!1]_q}(C_1q^{s\!+\!\frac{n+1}2}\!\!+\!
\!C_2q^{-s\!-\!\frac{n+1}2}\!+\!C_3)\!\!+\!\!\frac{\lambda_{2n+2}}{[2n\!+\!2]_q}C_3([n\!+\!1]_q\!-\!
[n\!+\!2]_q).\ea $$

Finally, we substitute the expression for $\tau_n(0)$ and use the
identities $$\ba{c} -[n+2]_q[n]_q-1+[n+1]_q[n+1]_q=0, \\[0.3cm]
-[n+2]_q(q^{n/2}+q^{-n/2})+k_q+[n+1]_q(q^{(n+1)/2}+q^{(n+1)/2})=0,
\ea$$ and the result follows. }

\bibliographystyle{plain}

\end{document}